\magnification\magstep1
\baselineskip 14pt
\def\biba{\par\parindent 20pt\hangindent 40pt}
\def\display#1:#2:#3\par{\par\hangindent #1 \noindent
			\hbox to #1{\hfill #2 \hskip .1em}\ignorespaces#3\par}
\def\disleft#1:#2:#3\par{\par\hangindent#1\noindent
			 \hbox to #1{#2 \hfill \hskip .1em}\ignorespaces#3\par}
\def\pfbox 
{{\ooalign{\hfil\lower.06ex 
 \hbox{$\scriptscriptstyle\smile$}\hfil\crcr
 \hfil\lower.7ex\hbox{\"{}}\hfil\crcr
 \mathhexbox20D}}}
\def\proof{\noindent{\sl Proof}.\enspace}

\centerline{\bf Stable Husbands}
\bigskip
\centerline{\sl Donald E. Knuth, Rajeev Motwani, and Boris Pittel}
\centerline{\sl  Computer Science Department, Stanford University}

\bigskip
{\narrower\smallskip\noindent
{\bf Abstract.}\enspace
Suppose $n$ boys and $n$ girls rank each other at random. We show that any
particular girl has at least $({1\over 2}-\epsilon)
\ln n$ and at most $(1+\epsilon)\ln n$
different husbands in the set of
all Gale/Shapley stable matchings defined by these rankings, with
probability approaching~1 as $n$, if $\epsilon$ is any positive 
constant. The proof emphasizes general methods that appear to be useful
for the analysis of many other combinatorial
algorithms.\footnote{}{%
This research 
was supported in part by the National Science Foundation under grant
CCR-86-10181,
and by Office of Naval Research contract
N00014-87-K-0502.}
\smallskip}

\bigskip
\noindent
{\bf 1. Introduction.}\enspace
This is a tale of $n$ girls and $n$ boys who play a game called ``stable
matching,'' invented by Gale and Shapley~[3].
Each player ranks each player of the opposite sex according to preference;
thus, there are $n$~permutations of the set of boys, representing the
preferences of the individual girls, and there are $n$~permutations
of the set of girls, representing the preferences of the individual boys.
The object of the game is for the boys and girls to match up so as to obtain
$n$~marriages that are {\it stable}, in the sense that no girl and boy prefer
each other to their current partners.

For example, suppose Alice, Brigitte, Cindy, and Debra play with Wilfred,
Xavier, Yuri, and Zeke; and suppose their preference rankings are as follows,
from favorite to least desired:
$$\vcenter{\halign{#\hfil\quad&$#$\hfil\qquad\qquad
&#\hfil\quad&$#$\hfil\cr
Alice likes&Y>X>Z>W&Wilfred likes&A>B>D>C\cr
Brigitte likes&X>W>Y>Z&Xavier likes&C>A>D>B\cr
Cindy likes&W>Y>X>Z&Yuri likes&B>D>A>C\cr
Debra likes&X>W>Z>Y&Zeke likes&B>A>C>D\cr
}}$$
The matching $(AW,BX,CY,DZ)$ is unstable, because for example
$A$ prefers~$Z$ to~$W$ and at the same time $Z$ prefers~$A$ to~$D$.
But the matching $(AZ,BW,CX,DY)$ is stable; most of the players 
are matched with a person other than their
first choice, but the objects of their affections don't want to change.
The given preferences also admit another stable matching, namely
$(AY,BW,CX,DZ)$. In this example only two matchings are stable.

The {\it stable husbands\/} of a girl are the boys she can be married to
in at least one stable matching. Thus, Alice's stable husbands in the
example are Yuri and Zeke. Brigitte has only one stable husband, namely
Wilfred; she likes Xavier better, but he can't stand her.

\vfill\eject

\noindent
{\bf 2. An algorithm.}\enspace
Gale and Shapley~[3]
gave a procedure to find a stable matching, given any set of preferences;
McVitie and Wilson~[10]
extended the method so that all stable matchings would be found. 
More recently, Gusfield~[5]
exploited the interesting lattice structure of stable matchings to
construct an elegant algorithm that simultaneously determines
the stable husbands of all girls in $O(n^2)$ steps. For our
purposes in the present paper
it suffices to consider a simplified variant of these
procedures, which finds the stable husbands of just one given girl~$G$.

The basic idea is to maintain partial matchings in which each boy who
currently has a partner is paired with his best possible choice among all
stable matchings of a certain class. 
One of the boys who doesn't have a current partner is temporarily called~$P$;
he will {\it propose\/} to one of the girls, and she will then decide whether
to accept or to reject his proposal (at~least for the time being).
The  role 
of~$P$ passes from boy to boy according to the following simple rules:

\medskip
\display 30pt:{\bf A0.}:
Initially all boys and girls are unpaired.

\display 30pt:{\bf A1.}:
If at least one boy has no current partner, let $P$ be one such boy.
Otherwise all boys and girls are already paired, and we have
a stable matching; output $G$'s partner~$S$ as one of her stable husbands,
then remove the pair~$GS$ from the current matching 
and let $P=S$. (Henceforth we will consider only stable matchings in
which $G$ is not married to~$S$.)

\display 30pt:{\bf A2.}:
If $P$ has already proposed to all the girls, terminate the algorithm. 
Otherwise let $H$ be the girl $P$~likes best among all those he hasn't
approached so far; $P$~now proposes to~$H$.

\display 30pt:{\bf A3.}:
If girl $H$ has already been proposed to by a boy she prefers to~$P$,
she rejects $P$'s~proposal. Otherwise, 
she accepts~$P$, and they become paired in the current matching;
her previous partner (if~any) now assumes the role of~$P$.
If she had no previous partner, the algorithm continues at step~A1,
otherwise it continues at~A2.\nobreak\quad\pfbox

\vfill\eject

For example, suppose we run this algorithm on the preference rankings given
in the introduction, 
always choosing the alphabetically least boy when there is a choice
in step~A1. 
Let the special girl~$G$ be Alice. Then the following events occur:

\smallskip
$$\vcenter{\halign{#\hfil\qquad&$#$\hfil\qquad
&$#\hfil$\qquad&$#$\hfil\qquad&#\hfil\cr
Step&\hbox{\rm Current matching}&P&H&Actions\cr
\noalign{\smallskip}
A1&&W\cr
A2&&W&A&$A$ accepts $W$\cr
A1&AW&X\cr
A2&AW&X&C&$C$ accepts $X$\cr
A1&AW,CX&Y\cr
A2&AW,CX&Y&B&$B$ accepts $Y$\cr
A1&AW,BY,CX&Z\cr
A2&AW,BY,CX&Z&B&$B$ rejects $Z$\cr
A2&AW,BY,CX&Z&A&$A$ accepts $Z$\cr
A2&AZ,BY,CX&W&B&$B$ accepts $W$\cr
A2&AZ,BW,CX&Y&D&$D$ accepts $Y$\cr
A1&AZ,BW,CX,DY&Z&&output $Z$\cr
A2&(AZ),BW,CX,DY&Z&C&$C$ rejects $Z$\cr
A2&(AZ),BW,CX,DY&Z&D&$D$ accepts $Z$\cr
A2&(AZ),BW,CX,DZ&Y&A&$A$ accepts $Y$\cr
A1&AY,BW,CX,DZ&Y&&output $Y$\cr
A2&(AY),BW,CX,DZ&Y&C&$C$ accepts $Y$\cr
A2&(AY),BW,CX,DZ&X&A&$A$ rejects $X$\cr
A2&(AY),BW,CY,DZ&X&D&$D$ accepts $X$\cr
A2&(AY),BW,CY,DX&Z&&terminate.\cr
}}$$

\smallskip
Each girl is paired with the boy who has made her the best offer so far,
except that the special girl Alice never has a partner in step~A2
after the first stable matching has been found. The notation `($AZ$)'
in this chart means that Alice has no current partner, but that her best
proposal so far has come from Zeke. In one place where `($AY$)' appears,
$A$~rejects $X$ even though she is currently unattached, because she
prefers~$Y$ to~$X$ and she will not lower her previous standards.

To prove that this algorithm finds all stable husbands of~$G$, let us suppose
for convenience that exactly one proposal is made per unit of time, so that
the $t\mskip1mu$th execution of step~A2 takes place at time~$t$. Denote by $Mt$
the set of all stable matchings such that $G$ is not married to any of
the $S$'s already output  before time~$t$. Then $M0$ is the set of all
stable matchings, and $tt'$ implies that $Mt\supseteq M{t'}$.
The correctness of the algorithm relies on the following crucial fact:
$$\eqalign{&\hbox{%
If girl $H$ rejects a suitor $R$ between time $t$ and $t+1$,}\cr 
\hbox{%
then }&\hbox{the pair $HR$
is not part of any stable matching in $M{t+1}$.}\cr}\eqno(\ast)$$
The proof is by induction on $t$. If~($\ast$) fails for the first
time at~$t$, suppose $H$ rejects~$R$ in step~A3 because she prefers~$Q$.
(Either $Q=P$ and $R$ is her previous best partner, or $R=P$ and $Q$ is
her previous best.) Then 
we are assuming that
$HR$ is part of a stable matching in $M{t+1}$,
and in this matching the stability condition tells us that boy~$Q$ must
be paired with some girl~$J$ he prefers to~$H$. But Quentin must then
 have proposed to Jane before he proposed to Helen, so he must have
been rejected by~$J$ at some time $t'<t$. Therefore, by~($\ast$) and induction,
$JQ$~is not part of a stable matching from~$M{t'}$, contradicting
the fact that $M{t'}\supseteq M{t+1}$. Rejection of~$R$ by~$H$ must
therefore have occurred not in step~A3 but in step~A1; in other words,
we must have $H=G$ and $R=S$, a stable husband. But then $M{t+1}$
does not include~$GS$, by definition. Therefore ($\ast$)~must be true.

The matchings found in step A1 must be stable. For if some girl~$H$
prefers boy~$U$ to her current partner, she has not yet been proposed
to by~$U$; hence he prefers his current mate. In fact, ($\ast$)~tells
us that the stable matchings found in~A1 are characterized by the
property that each boy has his best choice among all matchings in~$Mt$.
Thus, when the algorithm outputs~$S$, {\it each boy has his best choice among all
stable matchings such that\/ $G$ is paired with\/~$S$}.

The algorithm terminates when some boy has been rejected by all the girls.
According to~($\ast$), this happens at time~$t$ when $Mt=\emptyset$,
i.e., when all the stable husbands~$S$ of~$G$ have been output.

\vfill\eject

\noindent
{\bf 3. A random model.}\enspace
We wish to show that the algorithm just stated will produce at least $c\ln n$
outputs with probability approaching~1, if $c$ is any given constant
$<{1\over 2}$, assuming that the $n!^{2n}$ possible preference sequences
are selected uniformly at random.

The basic idea will be to use the principle of {\it late binding}, which also
has been called the principle of ``conservation of ignorance'' or ``deferred
decisions''---{\it le principe d'ajournement des d\'ecisions\/} in~[7].
Instead of fixing the preference sequences in advance, we simply let them
unfold to whatever extent the algorithm needs them as it runs. Thus,
whenever a boy is asked to propose, he proposes to a random girl chosen
uniformly from among the girls he hasn't tried yet. Whenever a girl receives
her $k$th~proposal, she accepts it with probability $1/k$. A~stochastic
process with these characteristics is equivalent to the original algorithm
running on random preference sequences, because it has the same
transition probabilities between states.

We can also simplify the algorithm further by assuming that each proposal
is uniformly random, as if each boy has ``amnesia''~[7]
and cannot remember any of the girls 
 he has previously asked. If it turns out that he has
just repeated himself, we will say that he has just made a {\it redundant
proposal\/}; such proposals are always rejected. The algorithm now reduces
to a fairly simple stochastic process, which uses the following data structures:

\medskip
\biba
$A1,\ldots,An=$ sets representing the girls proposed to so far
by boys~1 to~$n$.

\biba
$l=$ number of boys who have played the role of proposer.

\biba
$p=$ the boy who is currently proposing.

\biba
$h=$ the girl who is currently being proposed to.

\biba
$x1,\ldots,xn=$ boys who made the best offer so far to
girls 1 to~$n$, or zero if the girl has received no offer.

\biba
$k1,\ldots,kn=$ number of proposals received by girls 1 to~$n$.

\medskip
\display 30pt:{\bf B0.}:
Let $Aj=\emptyset$, $xj=0$, 
and $kj=0$, for $1jn$; also let $l=0$.

\display 30pt:{\bf B1.}:
If $l<n$, increase $l$ by~1 and let $p=l$.
Otherwise output~$xg$, where $g$ is the number of the special girl~$G$;
and let $p=xg$.

\display 30pt:{\bf B2.}:
Let $h$ be a random number, uniformly chosen between 1 and~$n$. (We say
that boy~$p$ has proposed to girl~$h$.) If $h\in Ap$ (i.e., if $p$'s
proposal is redundant), repeat this step. Otherwise replace $Ap$ by
$Ap\cup \{h\}$ and go on to step~B3.

\display 30pt:{\bf B3.}:
Increase $kh$ by one. With probability $1-1/kh$, return to~B2
(we say that girl~$h$ rejects the proposal). Otherwise interchange
$p\leftrightarrow xh$ (she accepts the proposal and her former
partner will have to propose to somebody else).
If the new value of~$p$ is zero, or if $h=g$ and at least one output
has already occurred, go back to step~B1; otherwise continue
with step~B2.\quad\pfbox

\medskip
Algorithm B faithfully models the previous Algorithm A on random input,
except for the redundant proposals. Notice that the new algorithm never
terminates; step~A2 stops when $P$ has nobody left to propose~to,
but step~B2 keeps making redundant proposals ad~infinitum when
$Ap=\{1,\ldots,n\}$. The details of Algorithm~B aren't extremely simple,
but we will see that certain aspects of its probable behavior are fairly
easy to analyze, in part because it never terminates.

\bigskip
\noindent
{\bf 4. Probabilistic preliminaries.}\enspace
Algorithm B can be regarded as a branching process, an infinite tree
with nodes at levels $t=1,2,3,\ldots$ corresponding to the $t\mskip1mu$th
time step~B2 is performed. Every node~$\alpha$ in this tree corresponds
to a unique path from the root, representing one of the possible behaviors
of the algorithm up to time~$t$. This path determines the values of the
data structures $(A1,\ldots,An,l,p1,h,x1,\ldots,xn,k1,\ldots,kn)$
at node~$\alpha$.

Every node $\alpha$ has $2n$ children $\alpha1^a,\alpha1^r,\ldots,\alphan^a,
\alphan^r$, where $\alphah^a$ and $\alphah^r$ represent the nodes following
a proposal that has been accepted or rejected by~$h$. If $h\in Ap$, the
transition probability from~$\alpha$ to~$\alphah^a$ is~0 and the 
transition probability from~$\alpha$ to~$\alphah^r$ is~$1/n$; this 
case corresponds
to a redundant proposal, which is always rejected. If $h\notin Ap$, the
transition probability from~$\alpha$ to~$\alphah^a$ is 
$1/(kh+1)n$
and the transition probability from~$\alpha$ to~$\alphah^r$ is
$kh/(kh+1)n$, where $kh$ is the data value that becomes $kh+1$ in
step~B3.

The probability ${\rm Pr}(\alpha)$ of node~$\alpha$ is the product of
the transition probabilities on the path from the root to~$\alpha$;
this is the probability that Algorithm~B will take the computational
path  represented by~$\alpha$. Since the transition probabilities from
each node to its children sum to~1, the sum of ${\rm Pr}(\alpha)$ for
all nodes~$\alpha$ on a given level~$t$ is~1.

We say that an event occurs at node $\alpha$ with {\it local probability\/}~$\rho$
if $\rho$ is the conditional probability of the event given that the
algorithm reaches~$\alpha$. Thus, for example, the local probability
that a proposal is accepted at~$\alpha$ is
$${1\over n}\sum{h\notin Ap}{1\over kh+1}\,,$$
the sum of the transition probabilities in which the event occurs. 

An event at $\alpha$
that depends only on the transition probabilities 
from~$\alpha$ to its children will be called
an {\it immediate event}.
 More general events may involve a sequence of node transitions in the subtree
below node~$\alpha$; all such events have local probabilities at~$\alpha$
as defined above. For example, we might speak of the local probability
that at most
five consecutive rejections immediately follow node~$\alpha$.
Local probabilities at~$\alpha$ are equivalent to unconditional 
probabilities in the branching process represented by the subtree whose
root is~$\alpha$.

Our proofs will often be based on a technique of probability estimation that
can conveniently be called the {\it principle of negligible perturbation}.
The idea will  be to change the transition probabilities between certain
nodes, obtaining a ``perturbed'' probability distribution Pr$'$ on which it is
relatively easy to compute the probability of some given event. Let $t$
be a level of the tree, and let $E$ be the set of all nodes at level~$t$
such that the given event is true. Let $C$ be the set of all nodes~$\alpha$
at level~$t$ whose probability has been perturbed somewhere along the path
from the root to~$\alpha$; thus, ${\rm Pr}(\alpha)={\rm Pr}'(\alpha)$
for all $\alpha\notin C$. Summing over all $\alpha\notin C$ and taking
complements tells us that ${\rm Pr}(C)={\rm Pr}'(C)$.
If ${\rm Pr}(C)$ is small, then the perturbation will have a negligible effect
on the probability of~$E$, because
$$\eqalign{|{\rm Pr}(E)-{\rm Pr}'(E)|
&=\left|\sum{\alpha\in E}{\rm Pr}(\alpha)-\sum{\alpha\in E}{\rm Pr}'
(\alpha)\right|\cr
\noalign{\smallskip}
&=\left|\sum{\alpha\in E\cap C}\bigl({\rm Pr}(\alpha)-{\rm Pr}'(\alpha)
\bigr)\right|\cr
\noalign{\smallskip}
&\sum{\alpha\in C}|{\rm Pr}(\alpha)-{\rm Pr}'(\alpha)|\cr
\noalign{\smallskip}
&\sum{\alpha\in C}|{\rm Pr}(\alpha)|
+\sum{\alpha\in C}|{\rm Pr}'(\alpha)|=2{\rm Pr}(C)\,.\cr}$$
Expected values can be estimated in a similar way.

(The principle of negligible perturbation seems almost absurdly simple, but we will
see that it simplifies our analyses in surprisingly nontrivial ways. The idea
is similar in spirit 
to Laplace's
method~[9] of asymptotic analysis, where integrals are estimated by
changing the integrand in unimportant portions of the domain. Another
kindred method is Wilkinson's
well-known technique of ``backward error analysis''~[13],
 in which numerical errors are conveniently studied by assuming that exact
answers have been obtained from approximate data; the actual situation,
in which approximate answers are calculated from exact data, is more difficult
to handle directly.)

Many of the proofs below are based on estimates of the tails of probability
distributions, using the following fundamental inequalities that we shall
call the {\it tail inequalities\/}: Let
$$P(z)=p0+p1z+p2z^2+\cdots =E(z^X)$$
be the probability generating function (pgf) for a random variable~$X$
that takes nonnegative integer values. Then
$$\eqalign{{\rm Pr}(Xr)\,&\,x^{-r}P(x)\qquad{\rm for}\enspace 0<x1\,;\cr
{\rm Pr}(Xr)\,&\,x^{-r}P(x)\qquad{\rm for}\enspace x1\,.\cr}$$
The proof is easy, since we have
$pkx^{-r}pkx^k$ when $0<x1$ and $kr$, and also
when $x1$ and $kr$.
In spite of this easy proof, the tail inequalities lead to
quite effective bounds because
we can often choose~$x$ to make $x^{-r}P(x)$ small.

(The history of these elementary inequalities takes us back to the early days of
probability theory. Bienaym\'e~[1] and Chebyshev~[12] observed that
$\Pr\bigl((X-\mu)^2r\bigr)E\bigl((X-\mu)^2\bigr)/r$ for all $r>0$.
Kolmogorov~[8] went further and remarked that $\Pr(Xr)E\bigl(f(X)\bigr)/s$ 
for any nonnegative function~$f(X)$, provided
that $E\bigl(f(X)\bigr)$ exists and $f(x)s>0$ for all $xr$. In
particular [8, equation 4.3.2], we get the second tail inequality when
$f(x)=e^{cx}$ and $c0$.
Chernoff~[2] pointed out the wide applicability of such estimates.)

\bigskip\noindent
{\bf 5. Probabilistic lemmas}. \enspace
Consider the behavior of Algorithm B as $n$. We will say that an event
occurs {\it almost surely}, or `a.s.', if the probability that it doesn't
happen is $o(1)$, i.e., if the probability of nonoccurrence approaches
zero as $n$. We will also say that an event occurs {\it quite surely},
or `q.s.', if the probability that it doesn't happen is superpolynomially
small, i.e., $O(n^{-K})$ for all fixed~$K$. If $p(n)$ is any polynomial
function, the sum of $O\bigl(p(n)\bigr)$ superpolynomially small probabilities
is superpolynomially small; hence if $m=O\bigl(p(n)\bigr)$ and if the events
$E1,\ldots,Em$ individually happen q.s., the combined event `$E1$ and
\dots\ and~$Em$' also happens q.s.

Let $N=\lfloor n^{1+\delta}\rfloor$, 
where $0<\delta<{1\over 2}$ is  a constant. Throughout
this section we shall consider only the first~$N$ proposals made by
Algorithm~B. Thus, probabilities of events are measured by summing
${\rm Pr}(\alpha)$ over all nodes~$\alpha$ at time $N+1$ such that the event
occurs as the algorithm follows the path to~$\alpha$.

\proclaim Lemma 1. Each girl q.s.\ receives at least ${1\over 2}n^{\delta}$
proposals and at most $2n^{\delta}$ proposals (including redundant ones).

\noindent
The statement of this lemma and those below is deliberately somewhat
ambiguous. One interpretation is that, if $g$ is any particular girl, she
q.s.\ receives the stated number of proposals. Another interpretation
is that q.s.\ all $n$ of the girls receive the stated number. The second
statement is a corollary of the first, because of the nature of `q.s.';
therefore we can prove each lemma using the first (weak) interpretation,
but we can apply each lemma by using the second (strong) interpretation.

\medskip
\proof
Let $g$ be one of the girls, and let $Ek$ be the event that the $k\/$th
proposal is to~$g$. This immediate event has local probability~${1\over n}$,
because each proposal in step~B2 is uniformly random. Therefore 
proposals to~$g$ are like Bernoulli trials with parameter~${1\over n}$, and
the pgf for
the total number of proposals received by~$g$ in the first $N$~levels
is simply
$$P(z)=\left({n-1+z\over n}\right)^N\,.$$
Let $r={1\over 2}n^{\delta}$. By the first tail inequality,
the probability
that $g$ receives at most $r$~proposals is at most
$$\left({1\over 2}\right)^{-r}P
\left({1\over 2}\right)=2^r
\left(1-{1\over 2n}\right)^{\lfloor 2nr\rfloor}2^r
\left(1-{1\over 2n}\right)^{2nr-1}2^{r+1}e^{-r}$$
since $1-xe^{-x}$, and this is superpolynomially small.

Similarly, if $r=2n^{\delta}$, the second tail inequality
tells us that $g$
receives~$r$ or more proposals with probability at most
$$2^{-r}P(2)=2^{-r}\left(1+{1\over n}\right)^{\lfloor{1\over 2}nr\rfloor}
2^{-r}\left(1+{1\over n}\right)^{{1\over 2}nr}
2^{-r}e^{{1\over 2}r}$$
since $1+xe^x$, again superpolynomially small.\quad\pfbox

\medskip
Let us say that a boy begins a {\it run of proposals\/} when he becomes
the proposer~$p$ in step~B1 or~B3; his run ends when one of his 
subsequent proposals is first accepted in step~B3.
In terms of the branching process, a~run continues when a transition is
from node~$\alpha$ to a ``rejected''
node of the form~$\alphah^r$, and it ends at
a transition from~$\alpha$ to an ``accepted'' node of the form~$\alphah^a$.

\proclaim
Lemma 2. Each boy q.s.\ begins at most $2n^{\delta}$ runs of proposals.

\proof
Let $b$ be one of the boys. His first run of proposals begins just
after $l$~increases to~$b$ in step~B1; his subsequent runs occur just after
$p$ is set to $xg=b$ in step~B1 or to $xh=b$ in step~B3. Thus, at most
two of his runs begin immediately after $p$ becomes~$b$ in step~B1.

The other runs occur when $p$ becomes $b=xh$ in step~B3; and this can happen
only if $h$ is the girl who accepted~$b$ at the end of his previous run.
Let $Et$ be the immediate event that the proposal at time~$t$ is to the girl
who has most recently accepted~$b$, or to girl~1 if $b$ has never yet been
accepted. Then the number of runs begun by~$b$ is at most 2~plus the number
of occurrences of~$Et$; in other words, $b$~can begin $r$ or more runs only
if $Et$ occurs $r-2$ or more times. But the local probability of~$Et$ 
is~${1\over n}$, so again we have the binomial pgf
$$P(x)=\left({n-1+x\over n}\right)^N$$
for the distribution of occurrences of~$Et$.

We now complete the proof as in Lemma 1, by setting $r=2n^{\delta}$;
the probability of~$r$ or more runs is at most $x^{2-r}P(x)$ for all $x>1$.
And we have seen that this bound is superpolynomially small when $x=2$.\quad
\pfbox

\proclaim Lemma 3. Each run q.s.\ contains at most $n^{\delta}(\log n)^2$
non-redundant proposals.

\proof
We will prove that for any fixed time $t$, $1tN$, a run starting at~$t$
q.s.\ has the stated property. Let $\alpha$ be any node at level~$t$,
and let $P(\alpha,m,t)$ be the local probability that the proposals
immediately following~$\alpha$ will include at least $m$~rejected 
non-redundant proposals before reaching time $N+1$ or before the first
acceptance, whichever comes first. Then we have the recursive formulas
$$P(\alpha,m,t)=\cases{1\,,&\hskip-40pt
if $m=0$;\cr
\noalign{\smallskip}
0\,,&\hskip-40pt
if $m>0$ and $t=N+1$;\cr
\noalign{\smallskip}
\displaystyle{%
\sum{h\in Ap}}{P(\alphah,m,t+1)\over n}+\displaystyle{\sum{h\notin Ap}}
{khP(\alphah^r,m-1,t+1)\over (kh+1)n}\,,&otherwise.\cr}$$

According to Lemma 1, we may assume that $kh2n^{\delta}$ for $1hn$.
(The validity of this assumption is discussed below.) Then it follows
by induction on $N+1-t$ that
$$P(\alpha,m,t)\left({2n^{\delta}\over 2n^{\delta}+1}\right)^m\,.$$
If we now choose $m=\lfloor n^{\delta}(\log n)^2\rfloor$, the local
probability that there are more than $m$~non-redundant proposals in a run
starting at~$\alpha$ is at most 
$$\exp\left({-m\over (2n^{\delta}+1)}\right)
=\exp\left(-{1\over 2}(\log n)^2+o(1)\right)\,.$$ 
Multiplying by ${\rm Pr}(\alpha)$
and summing over all~$\alpha$ on level~$t$ gives a total probability 
of at most
 $\exp\bigl(-{1\over 2}(\log n)^2+o(1)\bigr)$, which is
superpolynomially small.\quad\pfbox

\medskip
The previous proof uses a convenient simplification, indicated by the words
`According to Lemma~1, we may assume that \dots'. The assumption we are
making holds q.s., but it is not always true; moreover, it is a
probabilistic assertion about time $N+1$. So we should be careful that
we are not fallaciously using the future to influence probability
calculations in the past. 
A~rigorous justification can be made by
appealing to the principle of negligible perturbation: We simply recompute
the transition probabilities when the assumption $kh2n^{\delta}$ is
invalid. 

More precisely, if $\alpha$ is any node in the branching process, we let
the perturbed transition probabilities from~$\alpha$ to~$\alphah^a$ 
and~$\alphah^r$ be $1/(k'h+1)n$ and $k'h/(k'h+1)n$, 
respectively, where
$$k'h=\min(kh,2n^{\delta})\,.$$
The proof of Lemma~3 is valid for the perturbed branching process, using
$k'h$ in place of~$kh$ in the formula for $P(\alpha,m,t)$. Thus,
the proof establishes that each run in the perturbed branching process q.s.\
contains at most $n^{\delta}(\log n)^2$ non-redundant proposals. And
this same conclusion also holds q.s.\ in the unperturbed branching process,
because the probability of its falsity can increase by at most
$2{\rm Pr}(C)$, where $C$ is the condition that some transition probability
has been perturbed between the root and level $N+1$. Lemma~1 tells us
that ${\rm Pr}(C)$ is superpolynomially small, because the path to a node~$\alpha$
at level $N+1$ involves a perturbed transition probability only if some
girl in the state represented by~$\alpha$ has received more than~$2n^{\delta}$
proposals before time $N+1$.

\proclaim
Lemma 4. Each boy q.s.\ proposes to at most $2n^{2\delta}(\log n)^2$ girls.

\proof
Multiply the results of Lemmas 2 and 3.\quad\pfbox

\proclaim
Lemma 5. Each run q.s.\ contains at most $n^{\delta}(\log n)^2$ proposals.

\proof
Let $t$ be a fixed time, $1tN$, and let $\alpha$ be any node at level~$t$.
A~proposal is rejected with local probability 
$\sum{h\in Ap}1/n+\sum{h\notin Ap}kh/(kh+1)n$.

By the previous lemmas and the principle of negligible perturbation, we can
assume that $\|Ap\|2n^{2\delta}(\log n)^2$ and $kh2n^{\delta}$.
Let $\rho=\|Ap\|/n$.
Then the local probability of a run continuing one more step is at most
$$\rho+(1-\rho)\,{2n^{\delta}\over 2n^{\delta}+1}
={2n^{\delta}+\rho\over 2n^{\delta}+1}{2n^{\delta}+2n^{2\delta -1}(\log n)^2\over
2n^{\delta}+1}=\rho'\,.$$
(If the assumptions fail and the local probability is actually greater than
this number~$\rho'$,
 we can perturb it by artificially decreasing the probability of
rejection and increasing the probability of
acceptance. 
For example, we can define the transition probabilities from~$\alpha$
to~$\alphah^a$ and~$\alphah^r$ to be respectively $(1-\rho')/n$
and $\rho'/n$, for $1hn$.
The perturbed algorithm need not 
behave at all like the original algorithm does; 
 for example, a~boy's
redundant proposals might be accepted with positive probability.
The principle of negligible 
perturbations requires only that the 
nodes of the tree remain the same and that the
transition probabilities be consistent
with all assumptions of the proof.)

Since $\delta <{1\over 2}$, the local probability of $m$ consecutive redundant
or rejected proposals is at most
$$(\rho')^m<\left(1-{1\over 3n^{\delta}}\right)^m$$
for sufficiently large $n$. Hence we can complete
the proof as in Lemma~3.\quad\pfbox

\medskip
Notice that the principle of negligible perturbations has made it legitimate
for us, in this proof, to estimate probabilities of events that start
at time~$t$ by using assumptions that might fail at some future time~$>t$.
(Thus, $\|Ap\|$ might be $2n^{2\delta}(\log n)^2$ at the beginning of
a run but not at the end.) Arguments based on a weaker principle, which
would  require only that the assumptions hold at time~$t$, would be more
complicated; we would have to argue that $\|Ap\|$ cannot grow by more
than~1 at each time step, and our upper bound would be 
$\bigl(\rho'+m\big/\bigl(n(2n^{\delta}+1)\bigr)\bigr)^m$
instead of~$(\rho')^m$.

\proclaim
Lemma 6. Each boy q.s.\ makes at most $2n^{2\delta}(\log n)^2$ proposals.

\proof
Multiply the results of Lemmas 2 and 5.\quad\pfbox

\proclaim
Lemma 7. Each boy q.s.\ proposes to a given girl at most $\log n$ times.

\proof
Let $b$ be one of the boys and let $j$ be one of the girls. Perturb the
process so that after $b$ makes $n$~proposals, none of his subsequent
proposals 
has positive probability of being made
 to~$j$. This perturbation is negligible, because $b$ q.s.\
makes fewer than $n$~proposals (Lemma~6).

Furthermore, if $b$ hasn't made $n$ proposals by time $N+1$, pretend
that he continues proposing until he has done it $n$~times. This can only
increase the number of proposals he makes to~$j$.

The pgf for the total number of proposals by $b$ to~$j$ is then
$$P(z)=\left({n-1+z\over n}\right)^n\,,$$
because $b$ has amnesia; each of his $n$ proposals is uniform among the girls.
The probability that he has made more than $\log n$ of them to~$j$ is therefore
at most
$$(\ln n)^{-\log n}\left({n-1+\ln n\over n}\right)^n=\exp\bigl(-(\log n)
(\ln\ln n)+O(\ln n)\bigr)\,,$$
and this is superpolynomially small.\quad\pfbox

\smallskip
Incidentally, we have adopted here the convention of [4, p.~435] that
`log' is used for logarithms in contexts where the base is immaterial,
while `ln' denotes the special case of natural logs.

\proclaim
Lemma 8. Each girl q.s.\ receives at least ${1\over 2}n^{\delta}/\log n$
non-redundant proposals.

\proof
A girl receives q.s.\ ${1\over 2}n^{\delta}$ proposals by Lemma~1, but
at most $\log n$ from any one boy by Lemma~7.\quad\pfbox

\bigskip\noindent
{\bf 6. The main theorem.}\enspace
We are almost ready to show that Algorithm B a.s.\ produces $\Theta(\log n)$
outputs. (This final result will be ``almost sure'' but not ``quite sure.'')
But first we need to analyze the time of the first output, because steps~B1
and~B3 change their behavior at that time.

The first output occurs as soon as each of the $n$~girls has received at least
one proposal. We can prove that this q.s.\ happens long before 
time $N=\lfloor n^{1+\delta}\rfloor$:

\proclaim
Lemma 9. Let $N0=\lfloor n\ln n\ln\ln n\rfloor$. 
Each girl q.s.\ receives at least
one proposal and at most\/ $\ln n\,(\ln\ln n)^2$ proposals during the first
$N0$~steps.

\proof
The pgf for proposals to $g$ satisfies
$$P(x)=\left({n-1+x\over n}\right)^{N0}\exp\bigl((x-1)\ln n\ln\ln n
+o(1)\bigr)$$
for all real $x$. The probability that $g$ receives no proposal is
$P(0)\exp\bigl(-\ln n\ln\ln n+o(1)\bigr)$; the probability that she receives 
$\ln n\,(\ln\ln n)^2$ or more is at most
$$2^{-\ln n\,(\ln\ln n)^2}P(2)\exp\bigl(-(\ln 2)(\ln n)(\ln\ln n)^2
+\ln n\ln\ln n+o(1)\bigr)\,.$$
Both of these bounds are superpolynomially small.\quad\pfbox

\proclaim
Lemma 10. Let $\epsilon$ be a positive constant. A girl who has received
$m$~non-redundant proposals will accept at least\/ $(1-\epsilon)\ln m$
of them, with probability $1-O(m^{-\epsilon^2\!/2})$ as $m$.
She will accept at most\/ $(1+\epsilon)\ln m$ of them with probability
$1-O(m^{-\epsilon^2\!/2+\epsilon^3\!/6})$ as $m$.
And she will accept at most $m/(\ln m)^3$ of them with probability
$1-O\bigl(\exp\bigl(-m/(\ln m)^2+7m\,(\ln \ln m)/(\ln m)^3\bigr)\bigr)$
as $m$.

\proof
She accepts the $k\/$th with probability $1/k$, so the pgf for the total
number of acceptances is
$$P(z)=\left({z\over 1}\right)
\left({1+z\over 2}\right)
\ldots
\left({m-1+z\over m}\right)
={m-1+z\choose m}={1\over \Gamma(z)m^{\underline{1-z}}}\,.$$
(The notation $z^{\underline{w}}=z!/(z-w)!$ for factorial powers is
discussed in~[4, p.~211].)
The probability that she accepts fewer than $(1-\epsilon)\ln m$ is at most
$$(1-\epsilon)^{-(1-\epsilon)\ln m}P(1-\epsilon)={m^{-(1-\epsilon)\ln
(1-\epsilon)}\over \Gamma(1-\epsilon)\,m^{\underline{\epsilon}}}
{m^{\epsilon-\epsilon^2\!/2}\over \Gamma(1-\epsilon)\,
m^{\underline{\epsilon}}}\,,$$
and this is $O(m^{-\epsilon^2\!/2})$ because $m^{\underline{\epsilon}}=
m^{\epsilon}+O(m^{\epsilon-1})$. (See the answer to exercise 9.44 in~[4].)
Similarly, she accepts more than $(1+\epsilon)\ln n$ with probability at most
$$(1+\epsilon)^{-(1+\epsilon)\ln m}P(1+\epsilon){m^{-\epsilon-\epsilon^2\!/2
+\epsilon^3\!/6}\over\Gamma(1+\epsilon)\,m^{\underline{-\epsilon}}}
=O(m^{-\epsilon^2\!/2+\epsilon^3\!/6})\,.$$

The probability that she accepts more than $m0=m/(\ln m)^3$ of them is at most
$$\eqalign{m0^{-m0}P(m0)&=\exp\bigl(-m0\ln m0+\ln\Gamma(m+m0)
-\ln\Gamma(m0)-\ln m!\bigr)\cr
\noalign{\smallskip}
&=\exp\bigl(-m0(\ln m-6\ln\ln m+O(1))\bigr)\cr}$$
by Stirling's approximation.\quad\pfbox

\proclaim
Theorem. Assume that $n$ girls and $n$ boys have independent random preference
rankings, and let $G$ be one of the girls. Let $c$ be a constant $<{1\over 2}$
and let $C$ be a constant $>1$. Then $G$ a.s.\ has at least $c\ln n$ and at
most $C\ln n$ stable husbands.

\proof
The stable husbands of $G$ are output by Algorithm~A, which is equivalent
to Algorithm~B. The number of outputs is the number of times $g$ accepts
a proposal in Algorithm~B, minus the number of times she accepts a proposal
before the first output.

We have shown in Lemma 8 that $g$ will q.s. receive at least ${1\over 2}n^{\delta}/
\log n$ non-redundant proposals, among the first $n^{1+\delta}$ proposals made
by Algorithm~B, if $\delta$ is any constant between~0 and~${1\over 2}$. Therefore,
by the first estimate of Lemma~10, she will a.s.\ accept at least $(1-\epsilon)
\delta\ln n-O(\log\log n)$ proposals.

On the other hand, $g$ receives at most $n$ non-redundant proposals altogether.
Therefore, by the second estimate of Lemma~10, she will a.s.\ accept at most
$(1+\epsilon)\ln n$ of them.

Furthermore, by Lemma~9, the first output q.s.\ occurs before she has received
$m=\ln n\,(\ln\ln n)^2$ non-redundant proposals. Therefore (by the third estimate
of Lemma~10) she will accept at most
$${m\over(\ln m)^3}={\ln n\over \ln\ln n}\left(1+O\left({\log\log\log n\over
\log\log n}\right)\right)=o(\ln n)$$
before the first output, with probability
$$1-O\left(\exp\left(-\ln n+O\left({\log n\log\log\log n\over \log\log n}
\right)\right)\right)=1-{1\over n^{1+o(1)}}\,.$$

So the number of outputs will a.s.\ exceed $c\ln n$ 
and be less than~$C\ln n$
for all large~$n$, if we
choose $\delta$ and~$\epsilon$ so that $(1-\epsilon)\delta>c$ 
and $1+\epsilon <C$.
\quad\pfbox

\bigskip\noindent
{\bf 7. Remarks.}\enspace
Inspection of the proof of the theorem shows that the conclusion holds with
probability $1-O(n^{-\gamma})$, where $\gamma$ is any constant less than
both $(1-2c)^2\!/2$ and $(C-1)^2\!/2-(C-1)^3\!/6$.
We cannot improve this estimate to $1-O(n^{-1})$, because
there is probability $\sqrt{\,\ln n}/n$ that the first proposal to~$G$
will come from one of her $\sqrt{\,\ln n}$ favorite boys. In such a case
she can have at most $\sqrt{\,\ln n}$ stable husbands, because the first
 stable marriage found by Algorithm~A gives every girl her
{\it least\/} preferred stable husband.

Our theorem proves that random preferences a.s.\ guarantee an unbounded number
of stable matchings, since every stable husband is part of at least one
stable matching. Can it be shown that the a.s.\ lower bound of stable
matchings grows faster than this, say as $\Omega(\log n)^2$?
Pittel~[11]
has proved that the {\it expected\/} number of stable matchings is asymptotically
$e^{-1}n\ln n$. However, Pittel's theorem does not prove that a
large number of matchings will almost surely occur;
 constructions are known~[6]
where certain preference matrices give rise to at least $2^{n-1}$
 stable matchings, and such examples may be common enough to account
for the relatively high expected value.

\bigskip
\centerline{\bf References}

\medskip
\disleft 25pt:
[1]:
J. Bienaym\'e, ``Consid\'erations \`a l'appui de la d\'ecouverte de
Laplace sur la loi de probabilit\'e dans la m\'ethode des moindres
carr\'es,'' {\sl Comptes Rendus hebdomadaires des s\'eances de l'Acad\'emie des
 Sciences\/} (Paris) {\bf37} (1853), 309--324.

\smallskip
\disleft 25pt:
[2]:
Herman Chernoff, ``A measure of asymptotic efficiency for tests of a hypothesis
based on the sum of observations,''
{\sl Annals of Mathematical Statistics\/ \bf23} (1952), 493--507.

\smallskip
\disleft 25pt:
[3]:
D. Gale and L. S. Shapley, ``College admissions and the stability of
marriage,'' {\sl American Mathematical Monthly\/ \bf 69} (1962), 9--15.

\smallskip
\disleft 25pt:
[4]:
Ronald L. Graham, Donald E. Knuth, and Oren Patashnik, {\sl Concrete
Mathematics}, Addison-Wesley, 1989.

\smallskip
\disleft 25pt:
[5]:
Dan Gusfield, ``Three fast algorithms for four problems in stable marriage,''
{\sl SIAM Journal on Computing\/ \bf 16} (1987), 111--128.

\smallskip
\disleft 25pt:
[6]:
Robert W. Irving and Paul Leather, ``The complexity of counting stable
marriages,'' {\sl SIAM Journal on Computing\/ \bf 15} (1986), 655--667.

\smallskip
\disleft 25pt:
[7]:
Donald E. Knuth, 
{\sl Mariages Stables}, Les Presses de l'Universit\'e de Montr\'eal, 1976.

\smallskip
\disleft 25pt:
[8]:
A. Kolmogoroff, {\sl Grundbegriffe der Wahrscheinlichkeitsrechnung},
Springer, 1933. English translation by Nathan Morrison, {\sl Foundations
of the Theory of Probability}, Chelsea, 1950.

\smallskip
\disleft 25pt:
[9]:
P. S. La Place, ``M\'emoire sur les approximations des formules qui sont
fonctions de tr\`es grands nombres,'' {\sl M\'emoires de l'Academie royale
des Sciences de Paris\/} (1782), 1--88.
Reprinted in his {\sl \OE uvres Compl\`etes\/ \bf 10}, 207--291.


\smallskip
\disleft 25pt:
[10]:
D. G. McVitie and L. B. Wilson, ``The stable marriage problem,''
{\sl Communications of the ACM\/ \bf 14} (1971), 486--492.

\smallskip
\disleft 25pt:
[11]:
Boris Pittel, ``The average number of stable matchings,'' (submitted).

\smallskip
\disleft 25pt:
[12]:
P.-L. Tch\'ebyshef, 
``O srednikh velichinakh,'' {\sl Matematicheski\u\i\ Sbornik'\/ \bf 2}
(1867), 1--9; reprinted in his {\sl Polnoe Sobranie Sochineni\u\i},
volume~2, 431--437. French translation,
``Des valeurs moyennes,''
{\sl Journal de Math\'ema\-tiques pures et appliqu\'ees}, series~2, {\bf12}
(1867), 177--184; reprinted in his {\sl \OE uvres}, volume~1, 685--694.

\smallskip
\disleft 25pt:
[13]:
J. H. Wilkinson, {\sl Rounding Errors in Algebraic Processes}, Prentice-Hall,
1963.

\bye